\documentclass[12pt]{article}

\usepackage[utf8]{inputenc}
\usepackage[T2A]{fontenc}
\usepackage[english]{babel}

\usepackage{amsmath,amssymb,amsthm}
\usepackage{amsmath,amssymb,amsbsy,amsfonts,amsthm,latexsym,
            amsopn,amstext,amsxtra,euscript,amscd,mathrsfs,array,enumerate}
\usepackage{graphicx}
\usepackage[colorlinks=true, allcolors=blue]{hyperref}
\usepackage{tikz-cd}
\usepackage{csquotes}
\usepackage{mathtools}

\usepackage{float}

\newtheorem{theorem}{Theorem}
\newtheorem{proposition}{Proposition}
\newtheorem{lemma}[theorem]{Lemma}

\newtheorem{corollary}{Corollary}

\title{On the structure of solutions graph of \\ the Markoff--Hurwitz equation and recurrence sequences}

\date{}

\author{Ilya V. Vyugin\footnote{Department of Mathematics, HSE University, Usacheva Street 6, Moscow, 119048, Russia, e-mail: ilyavyugin@yandex.ru}}

\usepackage[a4paper,hmargin=2.5cm,vmargin=2.5cm]{geometry}

\begin{document}
 \maketitle

\begin{abstract} 
We prove that almost all solutions of the Markoff--Hurwitz equation over a residue field modulo $p$ can be obtained from one another by a chain of natural transformations. We also study recurrence sequences considered modulo prime $p$. 
\end{abstract}

\section{Introduction}

In 1879, Andrey Markoff studied in his master thesis in Saint Petersburg University the new Diophantine equation
\begin{equation}\label{eq:Markoff}
x^2+y^2+z^2=3xyz
\end{equation}
in integers $(x,y,z)\in \mathbb{Z}^3$. The natural solutions of this equation have found important applications in the theory of approximation of real numbers by rational ones, algebraic geometry, and other areas of mathematics. Markoff completely described all integer solutions of this equation. This equation is now called {\it Markoff equation}. Solutions $(x,y,z)$ to the Markoff equation are called {\it Markoff triples}. The elements of a Markoff triple are called {\it Markoff numbers}. Markoff triples and Markoff numbers have a number of properties. 
The open conjecture says that each natural Markoff number appears exactly once as the maximum number of a triple $(x,y,z)\in\mathbb{N}^3$ such that $x\geqslant y\geqslant z$ (see~\cite{Aigner}). The most important property of Markoff triples is that they form a tree-graph structure.


\bigskip

It is easy to see that if the triple $(x, y, z)\in \mathbb{Z}^3$ is a solution of the equation (\ref{eq:Markoff}), then the triples:
\begin{itemize}
\item[1.] all permutations of $(x,y,z)$;
\item[2.] $(-x,-y,z)$;
\item[3.] $(3yz-x,y,z)$
\end{itemize}
are also solutions of (\ref{eq:Markoff}). Moreover, if $(x,y,z)$ is a solution of (\ref{eq:Markoff}), then $(|x|,|y|,|z|)$ is also a solution. We can restrict ourselves to considering only natural triples $(x,y,z)\in \mathbb{N}^3$, since all integer solutions are either the zero solution $(0,0,0)$, or it is obtained from the natural solution by assigning two minuses, or it is natural. The natural solutions of equation (\ref{eq:Markoff}), considered up to permutations of components, form a tree graph with root $(1,1,1)$, where the solutions correspond to the vertices of the graph and the transformation $(x,y,z)\mapsto (3yz-x,y,z)$, which is an involution, corresponds to the edge leading from $(x,y,z)$ to $(3yz-x,y,z)$.

\medskip


Bourgain, Gamburd, and Sarnak \cite{BGS1} have recently started studying modulo $p$ reductions of the set of Markoff triples. Interesting problems arise if we consider the Markoff equation (\ref{eq:Markoff}) over a field $\mathbb{F}_p=\mathbb{Z}/p\mathbb{Z}$ of prime order $p$. It is easy to see that in this case transformations 1, 2 and 3 also take place and a Markoff graph of solutions is defined, which in this case is not a tree.

The Baragar; Bourgain, Gamburd, and Sarnak conjecture is as follows: {\it any solution of the Markoff equation over $\mathbb{F}_p$ can be obtained from two solutions $(0,0,0)$ and $(1,1,1)$ by a combination of the above transformations $1,2$ and $3$.}

William Chen gave a positive answer to this conjecture only in 2024 (see~\cite{Chen}) for sufficiently large prime $p$. The answer is based on two results obtained by methods that are very different from each other. 

In this case, it makes sense to consider a more general equation as well.
\begin{equation}\label{Markoff-ab}
x^2+y^2+z^2=Axyz+B,\qquad A\not=0,
\end{equation}
which can easily be reduced by substitution $(\tilde{x},\tilde{y},\tilde{z})=\frac{A}{3}(x,y,z)$ to the equation
\begin{equation}\label{Markoff-b}
x^2+y^2+z^2=3xyz+B',\quad B'=\frac{A^2}{9}B.
\end{equation}
For this equation, the transformations $1,2$ and $3$ hold (instead transformation 3 $(x,y,z)\mapsto (Ayz-x,y,z)$ is considered). The Markoff graph of solutions is defined, but does not contain triplets $(0,0,0)$, $(1,1,1)$ if $B\not= 0$.

For the Markoff graph of equations (\ref{Markoff-ab}),(\ref{Markoff-b}) Theorems \ref{th-BGS} and \ref{th-KSV} hold (see \cite{BGS1}, \cite{KSV} and Remark 3). We denote the Markoff graph of equation (\ref{Markoff-ab}) without the triple $(0,0,0)$ by $X^*(p)$. Let us denote connected components of the graph $X^*(p)$, by $C(p), D(p)$.

\begin{theorem}(Bourgain, Gamburd, Sarnak, 2016 \cite{BGS1})\label{th-BGS}
For any $\varepsilon>0$ and sufficiently large $p$, there exists a connected component $C(p)$ of the Markoff graph $X^{*}(p)$ such that
$$
|X^*(p)\setminus C(p)|\leqslant p^{\varepsilon}
$$
and for any connected component $D(p)$ of $X^{*}(p)$, we have
$$
|D(p)|>(\log p)^{1/3}.
$$
\end{theorem}
This theorem has been improved as follows:

\begin{theorem}(Konyagin, Shparlinskii, V'yugin 2022 \cite{KSV})\label{th-KSV}
There exists a connected component $C(p)$ of the generalized Markoff graph $X^{*}(p)$ such that
$$
|X^*(p)\setminus C(p)|\leqslant \exp((\log p)^{1/2+o(1)}),\quad p\to\infty
$$
and any connected component $D(p)$ satisfies the condition
$$
|D(p)|>c(\log p)^{7/9},
$$
where $c$ is the absolute constant.
\end{theorem}

Theorem 1.2.5 of Chen's paper \cite{Chen} shows that the number of elements of every non-zero connected component of $X^*(p)$ for equation (\ref{eq:Markoff}) is divisible by $p$. This, in combination with Theorem \ref{th-BGS} or \ref{th-KSV} resolves the conjecture of Bourgain, Gamburd, and Sarnak. We recall {\it giant component} by $C(p)$.
In this paper, we generalize the lower bound for the giant connected component to the Markoff--Hurwitz equation.

\medskip

The Markoff--Hurwitz equation 
\begin{eqnarray}\label{eq:MH}
x_1^2+\ldots+x_n^2=ax_1\ldots x_n
\end{eqnarray}
has basic solutions $(0,\ldots,0)$ and $(1,\ldots,1)$ (if $a=n$).
If the tuple $(x_1,\ldots,x_n)$ is a solution of equation (\ref{eq:MH}), then the tuples:
\begin{itemize}
\item[1.] all permutations of $(x_1,\ldots,x_n)$;
\item[2.] $(-x_1,-x_2,x_3,\ldots,x_n)$;
\item[3.] $(ax_2\ldots x_n-x_1,x_2,\ldots,x_n)$
\end{itemize}
are also solutions of (\ref{eq:MH}). These transformations define a graph of solutions to the Markoff--Hurwitz equation. The vertices of this graph correspond to solutions of equation (\ref{eq:MH}), the edges correspond to transformations 1,2 and 3. Let us denote the set of solutions $(x_1,\ldots,x_n)$ of equation (\ref{eq:MH}) by $\mathcal{H}$, and  the set of its solutions modulo $p$ by $\mathcal{H}(p)$.

Hurwitz proved that if $a=n$, a graph of non-zero solutions of the Markoff--Hurwitz equation in $\mathbb{Z}$ is connected. The question of connectivity of the solution graph of the Markoff--Hurwitz equation over a field $\mathbb{F}_p$ is still open. We are making progress in this problem in Theorem \ref{th-MH}.

\begin{theorem}\label{th-MH}
There exists a connected component $\hat{C}(p)$ of the Markoff-Hurwitz graph $H^{*}(p)$ such that
$$
|H^*(p)\setminus \hat{C}(p)|\leqslant \exp((\log p)^{1/2+o(1)})\leqslant p^{n-3+\varepsilon},\quad p\to\infty.
$$
\end{theorem}

\bigskip

In the second part of the paper we consider linear recurrence sequences modulo prime~$p$:
\begin{eqnarray}\label{recurr}
x_{n+2}=ax_{n+1}+bx_n,\quad x_1=x_1^0,\quad x_2=x_2^0,\quad x_n,a,b\in\mathbb{F}_p,\quad n\in\mathbb{N}.
\end{eqnarray}
These sequences are periodic. We denote the set of elements of the sequence (\ref{recurr}) by $C_{a,-b,x_1^0,x_2^0}$. We propose new upper bounds for the cardinalities of the intersections of the value sets of such sequences, their additive shifts, and multiplicative subgroups of the residue field. We apply these bounds to study the properties of modulo prime $p$ Fibonacci residues. 


The sequence $\{ x_n\}$ in the field $\mathbb{F}_p$ is a periodic sequence, since at some point two consecutive elements of the sequence will be repeated:
\begin{eqnarray*}
x_{L+1}=x_1,\quad x_{L+2}=x_2,
\end{eqnarray*}
$L$ is the length of the sequence period. We study intersections of subsets of $\mathbb{F}_p$ of the form:
$$
C_{a,-b,x_1^0,x_2^0};\quad C_{a,-b,x_1^0,x_2^0}+q=\{ x+q\mid x\in C_{a,-b,x_1^0,x_2^0}\};\quad G,
$$
where $G\subset\mathbb{F}_p$ is a subgroup and $b^k=1$, $a,b\in \mathbb{F}_p$, $k$ is a some natural number (see Prop.~\ref{prop1},\ref{prop2} and \ref{prop3}).

Estimates of this form have various applications. Estimates of intersections of recurrent sequences and their additive shifts can be used to estimate exponential sums. In addition, our estimates are applied to the study of properties of Fibonacci residue modulo prime sequences.

There are important applications of this kind of estimate, which underlies the proof of the well-known conjecture of Baragar, Bourgain, Gamburd, and Sarnak that the solutions of the Markoff equation over the field $\mathbb{F}_p$ are factorizations modulo $p$ of its integer solutions.

\section{Proof of the main result}

We need now to formulate some famous algebraic formulas. It is widely known that the number of points $N_p$ of an irreducible algebraic projective curve $\mathcal{C}$ of genus $g$ over the residue field $\mathbb{F}_p$ satisfies to Weil's bound
\begin{eqnarray}
|\#\mathcal{C}(\mathbb{F}_p)-(p+1)|\leqslant 2g\sqrt{p}.
\end{eqnarray}
It is easy to see that the number of points of an algebraic curve $\mathcal C=0\,:\, P(x,y)=0$ lying on the line at infinity does not exceed $\deg_x P+\deg_y P$.
It is well known that a curve that is given by an equation of degree two can have genus zero or one. 

\begin{lemma}\label{lem-irr}
Curve 
\begin{eqnarray}\label{curve-ab}
x_1^2+x_2^2-bx_1x_2+a=0
\end{eqnarray}
is irreducible if $a\not= 0$.
\end{lemma}

{\it Proof.} If the curve (\ref{curve-ab}) is reducible, then
$$
x_1^2+x_2^2-bx_1x_2+a=(x_1+\alpha x_2+\beta)(x_1+\gamma x_2+\delta)
$$
with some $\alpha,\beta,\gamma,\delta$. We see that $\alpha\delta=0$ and $\beta\gamma=0$, but $\beta\delta=a\not=0$, consequently, $\alpha=\gamma=0$ and $\alpha\gamma=1$ which is a contradiction. $\Box$

\bigskip

{\it Proof of Theorem \ref{th-MH}.} 
The idea of the proof is to consider three-dimensional reductions of the multidimensional Markoff--Hurwitz equation (\ref{eq:MH}). Three-dimensional reductions of the Markoff--Hurwitz equation have the form (\ref{Markoff-ab}) and Theorem \ref{th-KSV} (or \ref{th-BGS}) can be applied to them.

Consider sets of indices $\{ i_1,i_2,i_3\}$ and $\{ j_1,\ldots,j_{n-3}\}$
such that
$$
\{ i_1,i_2,i_3\}\cup \{ j_1,\ldots,j_{n-3}\}=\{ 1,2,\ldots,n\}.
$$
Let us denote $i=(i_1,i_2,i_3)$, $c=(c_{j_1},\ldots,c_{j_{n-3}})$. Recall the reduced Markoff--Hurwitz equation:
\begin{eqnarray}\label{rest-MH}
x_{i_1}^2+x_{i_2}^2+x_{i_3}^2+c_{j_1}^2+\ldots+c_{j_{n-3}}^2=ac_{j_1}\ldots c_{j_{n-3}}x_{i_1}x_{i_2}x_{i_3}
\end{eqnarray}
by $\mathcal{E}_{i,c}$. Equation (\ref{rest-MH}) is the equation (\ref{Markoff-ab}) with $A=ac_{j_1}\ldots c_{j_{n-3}}$, $B=-(c_{j_1}^2+\ldots+c_{j_{n-3}}^2)$. We apply Theorem \ref{th-KSV} (or \ref{th-BGS}) to equation (\ref{rest-MH}). The solutions graph of the reduced Markoff--Hurwitz equation $\mathcal{E}_{i,c}$ is a subgraph of the solutions graph $\mathcal{H}(p)$ of the Markoff--Hurwitz equation modulo $p$. We will prove that the giant components of Markoff graphs of all equations $\mathcal{E}_{i,c}$ for all $i$ and $c$ are connected to each other and form a giant connected component of the Markoff--Hurwitz graph $\mathcal{H}(p)$.

\begin{corollary}
For any reduced equation $\mathcal{E}_{i,c}$ with a non-zero set $(c_{j_1}\ldots c_{j_{n-3}})$ and for any $\varepsilon>0$ and sufficiently large $p$, there exists a giant connected component $C(p)$ of the Markoff graph $X^{*}(p)$ such that
$$
|X^*(p)\setminus C(p)|\leqslant  \exp((\log p)^{1/2+o(1)}) \leqslant p^{\varepsilon},\quad p\to \infty.
$$
\end{corollary}

This corollary is a direct consequence of Theorem \ref{th-BGS}.

\begin{lemma}\label{lem-giant-int}
Giant connected components of graphs of solutions two equations $\mathcal{E}_{i,c}$ and $\mathcal{E}_{i',c'}$:
\begin{eqnarray*}
\mathcal{E}_{i,c}\,:\, x_{i_1}^2+x_{i_2}^2+x_{i_3}^2+c_{j_1}^2+\ldots+c_{j_{n-3}}^2=ac_{j_1}\ldots c_{j_{n-3}}x_{i_1}x_{i_2}x_{i_3}
\end{eqnarray*}
and
\begin{eqnarray*}
\mathcal{E}_{i',c'}\,:\, x_{i_1}^2+x_{i_2}^2+x_{j_1}^2+c_{i_3}^2+c_{j_2}^2+\ldots+c_{j_{n-3}}^2=ac_{i_3}c_{j_2}\ldots c_{j_{n-3}}x_{i_1}x_{i_2}x_{j_1}.
\end{eqnarray*}
are intersected.
\end{lemma}

{\it Proof.} These two equations have common solutions and these solutions are solutions of the following equation:
\begin{eqnarray*}
x_{i_1}^2+x_{i_2}^2+c_{i_3}^2+c_{j_1}^2+\ldots+c_{j_{n-3}}^2=ac_{i_3}c_{j_1}\ldots c_{j_{n-3}}x_{i_1}x_{i_2}
\end{eqnarray*} 
in two variables $x_{i_1}$ and $x_{i_2}$.
If the tuple $(c_{i_3},c_{j_1},\ldots,c_{j_{n-3}})$ is non-zero, then this equation is irreducible (see Lemma \ref{lem-irr}). Thus, this equation has a number $N_c$ solutions such that:
$$
N_c>p-2\sqrt{p}-3.
$$
It follows from Weil's bounds.
We have in particular $N_c>p^{\varepsilon}$ ($\varepsilon\ll 1$). Tuples \\$(x_{i_1},x_{i_2},c_{i_1},c_{j_1},\ldots,c_{j_{n-3}})$ 
are also solutions of equations $\mathcal{E}_{i,c}$ and $\mathcal{E}_{i',c'}$ simultaneously. The number $N_c$ of solutions of this type is greater than $p^{\varepsilon}$ ($\varepsilon\ll 1$). This means that most of the equation solutions are solutions that belong to the giant connected components of equations $\mathcal{E}_{i,c}$ and $\mathcal{E}_{i',c'}$. This means that the giant connected components intersect, i.e. belong to the same connectivity component of the equation (\ref{eq:MH}). $\Box$

Applying Lemma \ref{lem-giant-int} first to the pair of reduced equations $\mathcal{E}_{i,c}$, $i=(i_1,i_2,i_3)$, $c=(c_{j_1},\ldots,c_{j_{n-3}})$ and
$\mathcal{E}_{i',c'}$, $i'=(i_1,i_2,j_1)$, $c'=(c_{i_1},c_{j_2},\ldots,c_{j_{n-3}})$. Secondly, we apply Lemma~\ref{lem-giant-int} to the pair of reduced equations $\mathcal{E}_{i',c'}$ and $\mathcal{E}_{i,c''}$, $i=(i_1,i_2,i_3)$, $c''=(c_{j_1}'',c_{j_2},\ldots,c_{j_{n-3}})$. Thus, in two steps, we can change the variable $c_{j_1}$ to any variable $c_{j_1}''$. So we can change any index $c_{j_l}$, $l=1,\ldots,n-3$. It is easy to see that applying Lemma~\ref{lem-giant-int} three times, the set of indices $i=(i_1,i_2,i_3)$ can be changed to any other set of indices $i'=(i_1',i_2',i_3')$. In all these transitions, the giant connectivity components become interconnected. We have shown that in $2(n-3)+3=2n-3$ steps (applications of Lemma~\ref{lem-giant-int}) one can pass from one reduced equation to any other reduced equation.

We have proved that the giant connected components of Markoff graphs of all reduced equations are connected to each other. Let us estimate from above the number of solutions of the Markoff--Hurwitz equation that are not included in this union of giant components. We have $C_n^3 p^{n-3}$ reduced equations. Each reduced equation has at most $p^{\varepsilon}$ (or $\exp((\log p)^{1/2+o(1)})$), $p\to\infty$ solutions that do not belong to the giant component. This means that no more than 
$$
C_n^3p^{n-3} \exp((\log p)^{1/2+o(1)})\quad  \text{or}\quad C_n^3 p^{n-3+\varepsilon},\quad p\to\infty,\quad C_n^3=\frac{n!}{(n-3)!3!}
$$
solutions can not be included in the union of all giant components. In this case, the total number of solutions to the Markoff--Hurwitz equation is within the limits:
$$
(p-3-2\sqrt{p})p^{n-2}<H^{*}(p)<(p+1+2\sqrt{p})p^{n-2}.
$$
$\Box$









\section{Intersections of recurrence sequences in $\mathbb{F}_p$}

Let us denote the set of values of the linear recurrence sequence:
\begin{eqnarray}\label{Rec-1}
x_{n+2}=a x_{n+1}-\varepsilon_{k} x_n,\quad x_1=x_1^0,\quad x_2=x_2^0,\quad n\in\mathbb{N},\quad \varepsilon_{k}^{k}=1
\end{eqnarray}
by $C_{a,\varepsilon_{k},x_1^0,x_2^0}$ and let us denote the set of values of the recurrence sequence:
\begin{eqnarray}\label{Rec-2}
y_{n+2}=b y_{n+1}-\varepsilon_{s} y_n,\quad y_1=y_1^0,\quad y_2=y_2^0,\quad n\in\mathbb{N},\quad \varepsilon_{s}^{s}=1
\end{eqnarray}
by $C_{b,\varepsilon_{s},y_1^0,y_2^0}$. Sequences (\ref{Rec-1}) and (\ref{Rec-2})
are considered in $\mathbb{F}_p$, where $p$ is a prime number.

The recurrence sequence (\ref{Rec-1}) is resolved by the formula:
\begin{eqnarray}\label{Rec-alpha-beta}
x_n=\alpha \lambda_1^n+\beta \lambda_2^{n},
\end{eqnarray}
where $\lambda_{1,2}=\frac{a\pm\sqrt{a^2-4\varepsilon_{k}}}{2}$ are the roots of the characteristic equation
\begin{eqnarray}\label{char-poly-x}
\lambda^2-a\lambda+\varepsilon_{k}=0.
\end{eqnarray}
Viet's formula gives us the relation $\lambda_1\lambda_2=\varepsilon_k$, i.e. $\lambda_2=\varepsilon_k \lambda_1^{-1}$. The constants $\alpha$ and $\beta$ are determined by the system:
\begin{eqnarray*}
\begin{cases}
\alpha \lambda_1+\beta \varepsilon_k \lambda_1^{-1}=x_1^0 \\
\alpha \lambda_1^2+\beta \varepsilon_k^2 \lambda_1^{-2}=x_2^0,
\end{cases}
\end{eqnarray*}
where $\alpha,\beta$ are constants. 
In our case, formula (\ref{Rec-alpha-beta}) can be rewritten as follows:
\begin{eqnarray*}
x_n=\alpha \lambda_1^n+\beta \varepsilon_k^n\lambda_1^{-n}.
\end{eqnarray*}
Let $n=lk+r$, then $x_{lk+r}=\alpha \lambda_1^{lk+r}+\beta \varepsilon_k^r\lambda_1^{-lk-r}=(\alpha \lambda_1^r)(\lambda_1^k)^l+ (\beta\lambda_1^{-r}\varepsilon_k^r)(\lambda_1^k)^{-l}$.

Divide the sequence $\{ x_n\}_{n=1}^{\infty}$ into $k$ subsequences, for each $r=0,\ldots,k-1$ we have the subsequence $C_{a,\varepsilon_{k},x_1^0,x_2^0}^r=\{ x_{lk+r}\}_{l=0}^{\infty}$. Powers of the element $\lambda_1$ generates a subgroup $L$ that lies in the quadratic extension of the field $\mathbb{F}_p$ extended by adjoining the roots of the characteristic polynomial (\ref{char-poly-x}), that is, $L\subset\mathbb{F}_{p^2}^*$ or $\mathbb{F}_p$ (in the case $\sqrt{a^2-4\varepsilon_k}\in\mathbb{F}_p$). Powers of the element $\lambda_1^k$ generates a subgroup $L$ of the order $\frac{|L|}{\text{gcd}(|L|,k)}$.

For the recurrence sequence (\ref{Rec-2}) the corresponding formulas are satisfied:
\begin{eqnarray}\label{Rec-gamma-delta}
y_n=\gamma \mu_1^n+\delta \mu_2^{n},
\end{eqnarray}
where $\mu_{1,2}=\frac{b\pm\sqrt{b^2-4\varepsilon_{s}}}{2}$ are the roots of the characteristic equation
\begin{eqnarray}\label{char-poly-y}
\mu^2-b\mu+\varepsilon_{s}=0
\end{eqnarray}
of the recurrence sequence (\ref{Rec-2}). Vieta's formulas give us the relation $\mu_1\mu_2=\varepsilon_s$, i.e. $\mu_2=\varepsilon_s \mu_1^{-1}$. The constants $\gamma$ and $\delta$ are determined by the system:
\begin{eqnarray*}
\begin{cases}
\gamma \mu_1+\delta \varepsilon_s \mu_1^{-1}=y_1^0 \\
\gamma \mu_1^2+\delta \varepsilon_s^2 \mu_1^{-2}=y_2^0.
\end{cases}
\end{eqnarray*}
The formula (\ref{Rec-alpha-beta}) can be rewritten as
\begin{eqnarray*}
y_n=\gamma \mu_1^n+\delta \varepsilon_s^n\mu_1^{-n}.
\end{eqnarray*}
Let $n=ls+r$, then $y_{ls+r}=\gamma \mu_1^{ls+r}+\delta \varepsilon_s^r\mu_1^{-ls-r}=(\gamma\mu_1^r)(\mu_1^s)^l+ (\delta\mu_1^{-r}\varepsilon_s^r)(\mu_1^s)^{-l}$.

The sequence $\{ y_n\}_{n=1}^{\infty}$ is also divided into $s$ subsequences, and for each $r=0,\ldots,s-1$ we have the subsequence $\{ y_{ls+r}\}_{l=0}^{\infty}$. The element $\mu_1$ generates a subgroup $H$ lying in the quadratic extension of the field $\mathbb{F}_p$ extended by adding roots of the characteristic polynomial (\ref{char-poly-y}), i.e. $H\subset\mathbb{F}_{p^2}^*$ or $\mathbb{F}_{p}^*$. The element $\mu_1^k$ generates a subgroup $H$ of the order $\frac{|H|}{\text{gcd}(|H|,s)}$.

We denote the sets of values of the subsequences $\{ x_{lk+r_1} \}_{l=1}^{\infty}$ and $\{ y_{ls+r_2} \}_{l=1}^{\infty}$ by $\tilde{C}_{a,\varepsilon_{k},x_1^0,x_2^0,r_1}$ and $\tilde{C}_{a,\varepsilon_{s},y_1^0,y_2^0,r_2}$, respectively.

\bigskip

The estimate of P. Corvaja and U. Zannier in~\cite{C-Z} and the corollary are uses in our proof.

\begin{theorem}(Corvaja, Zannier)\label{CZ-th}
Let $X$ be a smooth projective absolutely irreducible curve on a field $\kappa$ of characteristic $p$. Let $u,v\in \kappa(X)$ be rational functions that are multiplicatively independent modulo $\kappa^{*}$ with non-zero differentials, let $S$ be the set of their zeros and poles, let $\chi=|S|+2g-2$ be the Euler characteristic of $X\setminus S$. Then we have the estimate:
\begin{eqnarray}\label{poly-form}
\sum_{\nu\in X(\overline{\kappa})\setminus S}\min\{ \nu(1-u),\nu(1-v)\}\leqslant\max\left(3\sqrt[3]{2}(\deg u \deg v\chi)^{1/3},12\frac{\deg u \deg v}{p}\right),
\end{eqnarray}
where $\nu(f)$ denotes the multiplicity of $f$ vanishing at $\nu$.
\end{theorem}

From Corollary 2 \cite{C-Z} can be formulated as an estimate:
\begin{eqnarray}\label{curve-LH}
\#\{ (x,y)\mid (x,y)\in X, x\in L,y\in H \}\leqslant \max\left( 3\sqrt[3]{2}\chi^{1/3}|L|^{1/3}|H|^{1/3},12\frac{|L||H|}{p}\right),
\end{eqnarray}
which we shall use.

\begin{proposition}\label{prop1}
For the intersection of the sets of values of two sequences (\ref{Rec-1}) and (\ref{Rec-2}), the following inequality is satisfied:
$$
|C_{a,\varepsilon_{k},x_1^0,x_2^0}\cap C_{b,\varepsilon_{s},y_1^0,y_2^0}|\leqslant \max\left( 3\cdot 2^{2/3}(ks)^{2/3}|C_{a,\varepsilon_{k},x_1^0,x_2^0}|^{1/3} |C_{b,\varepsilon_{s},y_1^0,y_2^0}|^{1/3}, \frac{12}{p}|C_{a,\varepsilon_{k},x_1^0,x_2^0}| |C_{b,\varepsilon_{s},y_1^0,y_2^0}| \right).
$$
\end{proposition}

{\it Proof.}

The number of elements at the intersection $\tilde{C}_{a,\varepsilon_{k},x_1^0,x_2^0,r_1}\cap \tilde{C}_{b,\varepsilon_{k},y_1^0,y_2^0,r_2}$ of subsequences $\{ x_{lk+r_1} \}_{l=1}^{\infty}$ and $\{ y_{ls+r_2} \}_{l=1}^{\infty}$
does not exceed the number of pairs $(x,y)$ such that
$$
\alpha\lambda_1^{r_1} x+\frac{\beta\lambda_1^{-r_1}\varepsilon_k^{r_1}}{x}=\gamma\mu_1^{r_2} y+\frac{\delta\mu_1^{-r_2}\varepsilon_s^{r_2}}{y},
$$
that is, the numbers of points $(x,y)$ of the algebraic curve
\begin{eqnarray}\label{curveP12}
P_{r_1,r_2}(x,y)=\alpha\lambda_1^{r_1} x^2y+\beta\lambda_1^{-r_1}\varepsilon_k^{r_1}y-\gamma\mu_1^{r_2} xy^2-\delta\mu_1^{-r_2}\varepsilon_s^{r_2}x=0,
\end{eqnarray}
of genus $\leqslant 1$, such that $x\in L$, $y\in H$.
$$
C_{a,\varepsilon_{k},x_1^0,x_2^0}\cap C_{b,\varepsilon_{s},y_1^0,y_2^0}\subseteq\bigcup_{0\leqslant r_1<k,\, 0\leqslant r_2<s}\left(\tilde{C}_{a,\varepsilon_{k},x_1^0,x_2^0,r_1}\cap \tilde{C}_{b,\varepsilon_{k},y_1^0,y_2^0,r_2}\right).
$$
Then, respectively,
$$
|C_{a,\varepsilon_{k},x_1^0,x_2^0}\cap C_{b,\varepsilon_{s},y_1^0,y_2^0}|\leqslant \sum_{0\leqslant r_1<k,\, 0\leqslant r_2<s}|\tilde{C}_{a,\varepsilon_{k},x_1^0,x_2^0,r_1}\cap \tilde{C}_{b,\varepsilon_{k},y_1^0,y_2^0,r_2}|.
$$
Number of elements of each intersection $|\tilde{C}_{a,\varepsilon_{k},x_1^0,x_2^0,r_1}\cap \tilde{C}_{b,\varepsilon_{k},y_1^0,y_2^0,r_2}|$ is the number of points $(x,y)$ of the algebraic curve (\ref{curveP12}) such that $x\in L$, $y\in H$. Applying the estimate (\ref{curve-LH}), we obtain the estimate:

$$
|C_{a,\varepsilon_{k},x_1^0,x_2^0}\cap C_{b,\varepsilon_{s},y_1^0,y_2^0}|\leqslant ks\max\left(3\cdot 2^{2/3} |L|^{1/3} |H|^{1/3}, 12\frac{|L||H|}{p} \right)\leqslant
$$
$$
\leqslant\max\left( 3\cdot 2^{2/3}(ks)^{2/3}|C_{a,\varepsilon_{k},x_1^0,x_2^0}|^{1/3} |C_{b,\varepsilon_{s},y_1^0,y_2^0}|^{1/3}, \frac{12}{p}|C_{a,\varepsilon_{k},x_1^0,x_2^0}| |C_{b,\varepsilon_{s},y_1^0,y_2^0}| \right).
$$
The proof is complete.
$\Box$

\begin{proposition}\label{prop2}
For the intersection of the sets of values of two linear recurrence sequences (\ref{Rec-1}) and the subgroup $G\subset \mathbb{F}_p^*$ the inequality holds:
$$
|C_{a,\varepsilon_{k},x_1^0,x_2^0}\cap G|<
\max\left(3\sqrt[3]{4}k^{2/3}|C_{a,\varepsilon_{k},x_1^0,x_2^0}|^{1/3} |G|^{1/3},\frac{12}{p}|C_{a,\varepsilon_{k},x_1^0,x_2^0}| |G|\right).
$$
\end{proposition}

{\it Proof.} For the recurrence sequence (\ref{Rec-1}) the formulas from the previous theorem hold. The sequence $\{ x_n\}_{n=1}^{\infty}$ is split into subsequences $x_{lk+r}=(\gamma\lambda_1^r)(\lambda_1^k)^l+ (\delta\lambda_1^{-r}\varepsilon_k^r)(\lambda_1^k)^{-l}$,
$\lambda_{1}=\frac{a\pm\sqrt{a^2-4\varepsilon_{k}}}{2}$.

The number of elements of the intersection $\tilde{C}_{a,\varepsilon_{k},x_1^0,x_2^0,r}\cap G$ does not exceed the number of pairs $(x,y)$ such that
$$
\alpha\lambda_1^{r} x+\frac{\beta\lambda_1^{-r}\varepsilon_k^{r}}{x}=y, \quad x\in L,\quad y\in G.
$$
That is the number of points $(x,y)\in L\times G$ of the algebraic curve
$$
Q_{r}(x,y)=\alpha\lambda_1^{r} x^2+\beta\lambda_1^{-r}\varepsilon_k^{r} -xy=0
$$
of zero genus.
Then
$$
|C_{a,\varepsilon_{k},x_1^0,x_2^0}\cap G|\leqslant \sum_{0\leqslant r_1<k}|\tilde{C}_{a,\varepsilon_{k},x_1^0,x_2^0,r}\cap G|
\leqslant  k\max\left(3\sqrt[3]{4}|L|^{1/3} |G|^{1/3},\frac{12}{p}|L||G|\right)=
$$
$$
\leqslant\max\left(3\sqrt[3]{4}k^{2/3}|C_{a,\varepsilon_{k},x_1^0,x_2^0}|^{1/3} |G|^{1/3},\frac{12}{p}|C_{a,\varepsilon_{k},x_1^0,x_2^0}| |G|\right).
$$
The proof is complete.
$\Box$


\begin{proposition}\label{prop3}
For the intersection of the sets of values of the two linear recurrence sequences (\ref{Rec-1}) and the subgroup $G\subset \mathbb{F}_p^*$, the inequality holds:
$$
|C_{a,\varepsilon_{k},x_1^0,x_2^0}\cap (C_{a,\varepsilon_{k},x_1^0,x_2^0}-q) |\leqslant 
\max\left( 3\cdot 2^{2/3}k^{4/3}|C_{a,\varepsilon_{k},x_1^0,x_2^0}|^{2/3}, \frac{12}{p}|C_{a,\varepsilon_{k},x_1^0,x_2^0}|^2 \right).
$$
\end{proposition}
 
{\it Proof.} For the recurrence relation (\ref{Rec-1}) the formulas from the previous theorem hold. The sequence $\{ x_n\}_{n=1}^{\infty}$ is divided into subsequences $x_{lk+r}=(\gamma\lambda_1^r)(\lambda_1^k)^l+ (\delta\lambda_1^{-r}\varepsilon_k^r)(\lambda_1^k)^{-l}$,
$\lambda_{1}=\frac{a\pm\sqrt{a^2-4\varepsilon_{k}}}{2}$.

The number of elements of the intersection $\tilde{C}_{a,\varepsilon_{k},x_1^0,x_2^0,r_1}\cap (\tilde{C}_{a,\varepsilon_{k},x_1^0,x_2^0,r_2}-q)$ of subsequences $\{ x_{lk+r_1} \}_{l=1}^{\infty}$ and $\{ x_{lk+r_2}-q \}_{l=1}^{\infty}$ does not exceed the number of pairs $(x,y)$ such that
$$
\alpha\lambda_1^{r_1} x+\frac{\beta\lambda_1^{-r_1}\varepsilon_k^{r_1}}{x}=\alpha\lambda_1^{r_2} y+\frac{\beta\lambda_1^{-r_2}\varepsilon_k^{r_2}}{y}-q, \quad x\in G,\quad y\in G,
$$
that is, the numbers of points $(x,y)\in G\times G$ of the algebraic curve
$$
Q_{r_1,r_2}(x,y)=\alpha\lambda_1^{r_1} x^2y+\beta\lambda_1^{-r_1}\varepsilon_k^{r_1}y-\alpha\lambda_1^{r_2} xy^2-\beta\lambda_1^{-r_2}\varepsilon_k^{r_2}x+qxy=0.
$$
of genus $\leqslant 1$.
Then
$$
|C_{a,\varepsilon_{k},x_1^0,x_2^0}\cap (C_{a,\varepsilon_{k},x_1^0,x_2^0}-q)|\leqslant \sum_{0\leqslant r_1,r_2<k}|\tilde{C}_{a,\varepsilon_{k},x_1^0,x_2^0,r_1}\cap (\tilde{C}_{a,\varepsilon_{k},x_1^0,x_2^0,r_2}-q)|\leqslant
$$
$$
\leqslant
k^2\max\left(3\cdot 2^{2/3} |G|^{2/3},12\frac{|G|^2}{p} \right)\leqslant
\max\left( 3\cdot 2^{2/3}k^{4/3}|C_{a,\varepsilon_{k},x_1^0,x_2^0}|^{2/3}, \frac{12}{p}|C_{a,\varepsilon_{k},x_1^0,x_2^0}|^2 \right).
$$
The proof is complete.
$\Box$

\section{Applications to Fibonacci sequence}

An important special case of a linear recurrence relation is the sequence of Fibonacci numbers modulo integer $N$:
\begin{eqnarray}\label{Fibonacci}
f_{n+2}=f_{n+1}+f_n,\quad f_0=0,\quad f_1=1.
\end{eqnarray}
This sequence is periodic. The period of this sequence in $\mathbb{Z}/N\mathbb{Z}$ is called the Pisano period ($N$ is not necessarily prime). The length of the period is called the Pisano number and it is denoted by $\pi(N)$. There are known relations for $\pi(N)$, which almost completely express it in terms of the Pisano numbers of prime divisors of $N$:

\begin{itemize}

\item $\pi(mn)=\text{LCM} (\pi(m),\pi(n))$ for coprime $m$ and $n$;

\item $\pi(p^{k+1})=\pi(p^k)$ or $p\cdot\pi(p^k)$.

\end{itemize}

The Binet formula for Fibonacci numbers is as follows.
\begin{eqnarray}\label{Bine}
f_n=\frac{1}{\sqrt{5}}\left[\left(\frac{1+\sqrt{5}}{2}\right)^n+\left(-\frac{1+\sqrt{5}}{2}\right)^{-n} \right].
\end{eqnarray}
This formula holds for Fibonacci residues in $\mathbb{F}_p$, i.e. $x_n =\frac{1}{\sqrt{5}} \lambda^n + \frac{1}{\sqrt{5}}\frac{1}{(-\lambda)^n}$
where $\lambda=\frac{1+\sqrt{5}}{2}$ is the root of equation $\lambda^2-\lambda-1=0$. So, $\lambda\in\mathbb{F}_p$ or $\lambda\in\mathbb{F}_{p^2}=\{ l+k\sqrt{5}\mid k,l\in \mathbb{F}_p\}$.

Denote by $\mathcal{F}_p$ the set of Fibonacci residues modulo $p$. Then
$$
\mathcal{F}_p=\{ 5^{-1/2}[\lambda^{n}+(-\lambda)^{-n}]\mid n\in\mathbb{N} \}=\{ 5^{-1/2}[h+(-1)^nh^{-1}]\mid h\in H \},
$$
$\lambda^n$, $n\in\mathbb{N}$ generates some subgroup $H$ in $\mathbb{F}_{p^2}^*$, $(-\lambda)^n\in H\cup (-H)$.

\bigskip

In this paper, we only touch on periodicity modulo a prime number. The formulas for Pisano periods demonstrate the importance of the case of residues modulo a prime number.

\begin{corollary}
Let $G\in\mathbb{F}_p^*$ be subgroup. Then the inequality holds:
$$
|\mathcal{F}_p\cap G|<\max\left(6\sqrt[3]{2}(\pi(p))^{1/3} |G|^{1/3},\frac{12}{p}\pi(p)|G|\right).
$$
\end{corollary}

\begin{corollary}
For the sequence of Fibonacci residues (\ref{Fibonacci}) in $\mathbb{F}_p$, the inequality holds:
$$
|\mathcal{F}_p\cap (\mathcal{F}_p+q)|<
\max\left( 12 (\pi(p))^{2/3}, \frac{12}{p}(\pi(p))^2 \right).
$$
\end{corollary}

\end{document}